\documentclass[letterpaper, 10 pt, conference]{ieeeconf}  
\usepackage{lineno}
\usepackage[cmex10]{amsmath}
\usepackage{amscd,amsfonts,amsmath,amssymb,mathrsfs,pifont,stmaryrd,tipa}
\usepackage{array,cases,dsfont,graphicx,texdraw}
\usepackage{wrapfig}
\usepackage{graphics} 
\usepackage{epsfig} 
\usepackage{stfloats}
\usepackage{subfig}
\usepackage{hyperref}
\usepackage{color}

\newtheorem{Remark}{Remark}

\newenvironment{Proof}{\noindent{\em Proof:\/}}{\hfill $\Box$\par}
\newtheorem{Theorem}{Theorem}

\newtheorem{Assumption}{Assumption}

\input{tcilatex}                                                          

\IEEEoverridecommandlockouts                              
\title{\LARGE \bf
Distributed Nash Equilibrium Seeking for Games in Systems with Bounded Control Inputs}

\author{Maojiao Ye
\thanks{M. Ye is with the School of Automation, Nanjing University of Science and Technology, Nanjing 210094, P.R. China (Email: mjye@njust.edu.cn).}
\thanks{This work is supported by the National Natural Science Foundation of China (NSFC), No. 61803202, the Natural Science Foundation of
Jiangsu Province, No. BK20180455 and the Fundamental Research Funds for the Central Universities, No. 30920032203.}
}
\begin{document}

\maketitle
\thispagestyle{empty}
\pagestyle{empty}

\begin{abstract}
Noticing that actuator limitations are ubiquitous in practical engineering systems, this paper considers Nash equilibrium seeking for games in systems where the control inputs are bounded. More specifically, first-order integrator-type systems with bounded control inputs are firstly considered and two saturated control strategies are designed to seek the Nash equilibrium of the game. Then, second-order integrator-type systems are further considered. In this case, a centralized seeking strategy is firstly proposed without considering the boundedness of the control inputs, followed by a distributed counterpart. By further adapting a saturation function into the distributed Nash equilibrium seeking strategy, the boundedness of the control input is addressed. In the proposed distributed strategies, consensus protocols are included for information sharing and the saturation functions are utilized to construct bounded control inputs. The convergence results are analytically studied by Lyapunov stability analysis. Lastly, by considering the connectivity control of mobile sensor networks, the proposed methods are numerically verified.
\end{abstract}


\section{Introduction}

Games are attracting growing interests from researchers in the multi-agent communities for the analysis of multi-agent systems in recent years \cite{Parsons02}. For example, consensus was accomplished by utilizing cooperative game theory in \cite{Semsar-Kazerooni09}. Differential games were applied to solve distributed optimal tracking control of multi-agent systems with external disturbance in \cite{JiaoAT16}. The works in  \cite{MardenTSNC09}-\cite{YIJCD} linked games to cooperative control and optimization of multi-agent systems, respectively. In \cite{MaIJRNC18}, the consensus analysis for a class of hybrid multi-agent systems was conducted based on a noncooperative game. These works motivate us to take the physical constraints of multi-agent systems into consideration for Nash equilibrium seeking problems. The concerned constraints include but are not limited to communication issues, input saturation, system dynamics and action constraints. Recent years witnessed the trials made by researchers to accommodate system dynamics (see, e.g., \cite{StankovicTAC12}\cite{IbrahimCDC18}\cite{DengAT19}\cite{FilippoECC19}), the communication issues for games in distributed networks (see, e.g., \cite{YETcyber18}\cite{DengTNNLS}) and action constraints (see, e.g., \cite{DengTNNLS}\cite{LuTcyber}). For example, following the ideas presented in \cite{YETcyber17}-\cite{YETAC19} to establish distributed Nash equilibrium seeking strategies by utilizing consensus algorithms and the gradient search, communication constraints were accommodated in \cite{YETcyber18}. Moreover, weight-balanced digraphs were considered in \cite{DengTNNLS}. Games in linear systems and Euler-Lagrange systems  were considered in \cite{StankovicTAC12} and \cite{DengAT19}, respectively. Un-modeled dynamics and disturbances were addressed in \cite{YeTCYdoi}. In addition,  \cite{IbrahimCDC18} and \cite{FilippoECC19} focused on second-order dynamics. High-order games were considered in \cite{RomanoTCNS20}\cite{RomanoECC19}, where internal-model-based seeking strategies were proposed to achieve distributed Nash equilibrium seeking.
Generalized Nash equilibrium seeking, which concerns with action constraints among the players, was handled in \cite{Auto16}\cite{TatarenkoTCA}\cite{Belgioioso17}\cite{PersisTAC20}. Besides, the extremum seeking based perspectives in \cite{FrihaufTAC12}, the gossip algorithms in \cite{KoshalOR16}\cite{SalehisadaghianiAT16}, the passivity perspectives in \cite{Gadjov19} and the integral dynamics in \cite{PersisTAC20}\cite{PersisAT19} also provided insightful ideas to achieve Nash equilibrium seeking. However, actuation limitations are not considered in these works.

As many engineering systems are subject to actuator limitations (e.g., robotic manipulators \cite{HeTSMC16}, spacecraft \cite{Boskovic}, hard disk drive servo systems \cite{ChenTAC03}, just to name a few), the boundedness of control inputs appears to be a problem that is both practically and theoretically concerned. The study for systems with bounded control inputs has a rich history.  For example, input-saturated linear systems were considered in \cite{KapporAT98} based on an anti-windup design. Backstepping approaches were employed for developing robust adaptive control strategies to accommodate uncertain nonlinear systems subject to input saturation \cite{WenTAC11}. Two-player zero-sum games with non-quadratic payoffs were employed to solve the $H_{\infty}$ control of systems with bounded control inputs in \cite{Modares14}. Moreover, with the development of multi-agent systems, consensus problems in input-saturated multi-agent systems have attracted a lot of attention. The authors in \cite{Mengscl13} dealt with leader-following consensus of linear multi-agent systems with input saturation. Global consensus of saturated discrete-time systems was addressed in \cite{YangAT14}. Optimal consensus for multi-agent systems with bounded control inputs was investigated in \cite{XieSCL17}-\cite{QiuIJRNC}. However, Nash equilibrium seeking for games in systems with bounded controls has not been addressed yet, though it is a problem of great interest.


Inspired by the above observations, we intend to design Nash equilibrium seeking strategies for games in both first-order and second-order integrator-type systems in which the controls are bounded. The considered problem is challenging as the saturation function would introduce high nonlinearity into the closed-loop system. Moreover, the nonlinearity would result in difficulties on the design of the Nash equilibrium seeking algorithms, the establishment of the Lyapunov functions and the corresponding stability analyses.
In summary, with part of the manuscript presented in \cite{YEICCA19}, this paper contributes in the following aspects: 1). Distributed Nash equilibrium seeking for games in systems with bounded control inputs is considered in this paper. First-order integrator-type systems are firstly considered, in which both the saturated gradient play and a distributed strategy are investigated. Then, second-order integrator-type systems are explored. A centralized algorithm is firstly proposed without considering the boundedness of the control inputs, followed by two distributed seeking schemes. 2). The convergence results of the proposed Nash equilibrium seeking strategies are analytically investigated. It is proven that the proposed seeking strategies would enable the players' actions to asymptotically converge to the Nash equilibrium under the given conditions.

\textbf{Notations:} In the remainder, we use $\mathbb{R}$ to denote the set of real numbers. The notation $[h_i]_{vec}$ is defined as $[h_{i}]_{vec}=[h_1,h_2,\cdots,h_N]^T$ and $diag\{h_{ij}\}$($diag\{h_i\}$) for $i,j\in\{1,2,\cdots,N\}$ denotes a diagonal matrix whose diagonal elements are $h_{11},h_{12},\cdots,h_{1N},h_{21},\cdots,h_{NN},$ ($h_1,h_2,\cdots,h_N$), successively. For a symmetric matrix $Q\in \mathbb{R}^{N\times N},$ $\lambda_{min}(Q)$ denotes the minimum eigenvalue of $Q$. Moreover, $\otimes$ is the Kronecker product. The notation $\text{min}\{a,b\}=a$ if $a\leq b$ and $\text{min}\{a,b\}=b$ if $a>b.$ In addition, $\tilde{H}=[h_{ij}]$ defines a matrix whose $(i,j)$th entry is $h_{ij}.$


\section{Problem Formulation}\label{p1_res}
Consider a game with $N$ players whose dynamics are governed by
\begin{equation}
x_i^n=u_i,
\end{equation}
where $x_i\in \mathbb{R}$ is the action of player $i$ and $u_i\in\mathbb{R}$ is the control input that satisfies  $|u_i|\leq \bar{U}.$ Moreover, $x_i^n$ denotes the $n$th-order time derivative of $x_i$ and in the subsequent section, $n=1$ and $n=2$ will be investigated successively. Let $f_i(\mathbf{x}),$ where $\mathbf{x}=[x_1,x_2,\cdots,x_N]^T,$ be the cost function of player $i$ and $\{1,2,\cdots,N\}$ denotes the set of $N$ players. This paper aims to design the bounded controls to seek the Nash equilibrium $\mathbf{x}^*=(x_i^*,\mathbf{x}_{-i}^*)$ on which
\begin{equation}
f_i(x_i^*,\mathbf{x}_{-i}^*)\leq f_i(x_i,\mathbf{x}_{-i}^*),
\end{equation}
for $x_i\in \mathbb{R},i\in\{1,2,\cdots,N\}$ and $\mathbf{x}_{-i}=[x_1,x_2,\cdots,x_{i-1},x_{i+1},\cdots,x_{N}]^T$.

The following conditions will be utilized to establish the convergence results.
\begin{Assumption}\label{Assu_1}
The players' cost functions are $\mathcal{C}^2$ functions.
\end{Assumption}
\begin{Assumption}\label{a3}
The players are equipped with a communication graph $\mathcal{G},$ which is undirected and connected.
\end{Assumption}

\begin{Assumption}\label{Assu_2}
\cite{YETAC17}\cite{YeAT18} There exists a positive constant $m$ such that
\begin{equation}\label{cond_1}
(\mathbf{x}-\mathbf{z})^T(\bar{\mathcal{P}}(\mathbf{x})-\bar{\mathcal{P}}(\mathbf{z}))\geq m||\mathbf{x}-\mathbf{z}||^2,
\end{equation}
for all $\mathbf{x},\mathbf{z}\in \mathbb{R}^N$. Note that in \eqref{cond_1}, $\bar{\mathcal{P}}(\mathbf{x})=\left[\nabla_if_i(\mathbf{x})\right]_{vec}$ and $\nabla_if_i(\mathbf{x})=\frac{\partial f_i(\mathbf{x})}{\partial x_i}$.
\end{Assumption}

\begin{Assumption}\label{Ass4}
The elements in $H(\mathbf{x}),$ defined as $H(\mathbf{x})=\left[\frac{\partial^2 f_i(\mathbf{x})}{\partial x_i\partial x_j}\right],$ are bounded for $\mathbf{x}\in \mathbb{R}^N.$
\end{Assumption}

\begin{Remark}
From Assumption \ref{Assu_2}, it can be obtained that for each fixed $\mathbf{x}_{-i},$ $f_i(x_i,\mathbf{x}_{-i})$ is strongly convex in $x_i$ and $H^T(\mathbf{x})+H(\mathbf{x})\geq 2m I$
for $\mathbf{x}\in \mathbb{R}^N$ by Proposition 2.3.2 in \cite{Facchinei}. Moreover, under Assumption \ref{Assu_2}, the game admits a unique Nash equilibrium by Theorem 2.3.3 in \cite{Facchinei} and the players' actions are at the Nash equilibrium if and only if $\bar{\mathcal{P}}(\mathbf{x})=\mathbf{0}_N$ \cite{YETAC17}.  Assumption \ref{Ass4} indicates that for each $i\in\{1,2,\cdots,N\},$ $\nabla_if_i(\mathbf{x})$ is globally Lipschitz. Moreover, it's worth noting that Assumption \ref{Ass4} is utilized for the development of global convergence results, and without this condition, weaker convergence results can be obtained.
\end{Remark}
\begin{Remark}
Note that compared with our previous works in \cite{YETcyber18}\cite{YETcyber17}-\cite{YETAC19}, it is required that $|u_i|\leq \bar{U}$ in this paper. Due to the high nonlinearity introduced by the boundedness of controls, the establishments of the seeking strategies and the associated  Lyapunov stability analysis would be challenging. Moreover, the graph related definitions utilized in the paper follow those in \cite{YETAC17} and are omitted directly in this paper due to space limitation.
\end{Remark}

\section{Main Results}\label{min_re}

In this section, Nash equilibrium seeking for games in which the players are of first-order integrator-type dynamics and second-order integrator-type dynamics will be successively investigated.

\subsection{First-order integrator-type systems}
In this section, we consider games in which the players' actions are governed by
\begin{equation}\label{first}
\dot{x}_i=u_i,i\in\{1,2,\cdots,N\}.
\end{equation}
In the following, saturated gradient play will be firstly considered, followed by a distributed seeking strategy.
\subsubsection{Saturated gradient play}
To seek the Nash equilibrium of the game, we suppose that the players update their actions according to
\begin{equation}\label{satur_gra}
\dot{x}_i=-\rho_{\bar{U}}\left(\nabla_if_i(\mathbf{x})\right),
\end{equation}
where $i\in\{1,2,\cdots,N\},$ and  $\rho_{\bar{U}}(\eta_i)=sgn(\eta_i)\min\{|\eta_i|,\bar{U}\}.$

\begin{Theorem}\label{Lem_1}
The Nash equilibrium of the game is globally asymptotically stable under \eqref{satur_gra} given that Assumptions \ref{Assu_1} and \ref{Assu_2} are satisfied.
\end{Theorem}
\begin{Proof}
Let
$V\left(\bar{\mathcal{P}}(\mathbf{x})\right)=\sum_{i=1}^N \int_{0}^{\nabla_if_i(\mathbf{x})} \rho_{\bar{U}}(t)dt$
be the Lyapunov candidate function.
Then, if $0\leq \nabla_if_i(\mathbf{x})<\bar{U}$, $\int_{0}^{\nabla_if_i(\mathbf{x})}\rho_{\bar{U}}(t)dt=\frac{1}{2}\left(\nabla_if_i(\mathbf{x})\right)^2,$
and if $\nabla_if_i(\mathbf{x})\geq \bar{U},$ $\int_{0}^{\nabla_if_i(\mathbf{x})}\rho_{\bar{U}}(t)dt=\frac{1}{2}\bar{U}^2+\left(\nabla_if_i(\mathbf{x})-\bar{U}\right)\bar{U}.$
Therefore, $\int_{0}^{\nabla_if_i(\mathbf{x})}\rho_{\bar{U}}(t)dt>0$ for $\nabla_if_i(\mathbf{x})>0$ and  $\int_{0}^{\nabla_if_i(\mathbf{x})}\rho_{\bar{U}}(t)dt\rightarrow +\infty$ as $\nabla_if_i(\mathbf{x})\rightarrow +\infty.$ In addition, if $\nabla_if_i(\mathbf{x})<0,$ $\int_{0}^{\nabla_if_i(\mathbf{x})}\rho_{\bar{U}}(t)dt=-\int_{\nabla_if_i(\mathbf{x})}^{0}\rho_{\bar{U}}(t)dt
=\int_{0}^{\left|\nabla_if_i(\mathbf{x})\right|}\rho_{\bar{U}}(t)dt.$
Hence, $\int_{0}^{\nabla_if_i(\mathbf{x})}\rho_{\bar{U}}(t)dt>0$ as well for $\nabla_if_i(\mathbf{x})<0,$ and $\int_{0}^{\nabla_if_i(\mathbf{x})}\rho_{\bar{U}}(t)dt\rightarrow +\infty,$ as $\nabla_if_i(\mathbf{x})\rightarrow -\infty.$ Moreover, if  $\nabla_if_i(\mathbf{x})=0,$ it is clear that $\int_{0}^{\nabla_if_i(\mathbf{x})}\rho_{\bar{U}}(t)dt=0.$ Recalling that $\bar{\mathcal{P}}(\mathbf{x})=\left[\nabla_if_i(\mathbf{x})\right]_{vec},$ it can be derived that the Lyapunov candidate function is positive definite with respect to $\bar{\mathcal{P}}(\mathbf{x}).$ Moreover, if $||\bar{\mathcal{P}}(\mathbf{x})||\rightarrow +\infty,$ there exists at least one player $j$ whose gradient value satisfies $\left|\nabla_jf_j(\mathbf{x})\right|\rightarrow +\infty,$ and  hence, $\int_{0}^{\nabla_jf_j(\mathbf{x})}\rho_{\bar{U}}(t)dt\rightarrow +\infty.$ As for each $i\in\{1,2,\cdots,N\},$ $\int_{0}^{\nabla_if_i(\mathbf{x})}\rho_{\bar{U}}(t)dt\geq 0,$ we can conclude that $V(\bar{\mathcal{P}}(\mathbf{x}))\rightarrow +\infty$ as $||\bar{\mathcal{P}}(\mathbf{x})||\rightarrow +\infty,$ i.e., the Lyapunov candidate function is radially unbounded with respect to $\bar{\mathcal{P}}(\mathbf{x})$.

Taking the time derivative of $V$ gives $\dot{V}=\sum_{i=1}^N \rho_{\bar{U}}\left(\nabla_if_i(\mathbf{x})\right) \left(\frac{\partial}{\partial \mathbf{x}}\left(\nabla_if_i(\mathbf{x})\right)\right)^T \dot{\mathbf{x}}
=-\left[\rho_{\bar{U}}\left(\nabla_if_i(\mathbf{x})\right)\right]_{vec}^T H(\mathbf{x})\left[\rho_{\bar{U}}\left(\nabla_if_i(\mathbf{x})\right)\right]_{vec}.$
By Assumption \ref{Assu_2}, $H^T(\mathbf{x})+H(\mathbf{x})\geq 2mI$. Therefore, $\dot{V}\leq - m\left|\left| \left[\rho_{\bar{U}}\left(\nabla_if_i(\mathbf{x})\right)\right]_{vec}\right|\right|^2.$
Hence, $\left|\left| \rho_{\bar{U}}\left(\nabla_if_i(\mathbf{x})\right)\right|\right|\rightarrow 0$ for all $i\in\{1,2,\cdots,N\}$ as $t\rightarrow +\infty$. Noticing that by Assumption \ref{Assu_2}, $\bar{\mathcal{P}}(\mathbf{x})=\mathbf{0}_N$ if and only if $\mathbf{x}=\mathbf{x}^*$, we can conclude that $||\mathbf{x}-\mathbf{x}^*||\rightarrow 0$ as $t\rightarrow +\infty.$
\end{Proof}

In Theorem \ref{Lem_1}, the convergence property of the saturated gradient play in \eqref{satur_gra} is investigated. However, the saturated gradient play is not suitable for distributed games as all the players' actions are contained in the gradient information.  Therefore, we further investigate the Nash equilibrium seeking problem under distributed networks in the subsequent section.
\subsubsection{Consensus-based distributed Nash equilibrium seeking}
To achieve Nash equilibrium seeking in distributed networks, we suppose that the players can communicate with each other via communication graph $\mathcal{G}$. Then, the Nash equilibrium seeking strategy can be designed as
\begin{equation}\label{bound_dis}
\begin{aligned}
\dot{x}_i=&-\rho_{\bar{U}}\left(\nabla_if_i(\mathbf{y}_i) \right),\\
\dot{y}_{ij}=&-\theta_{ij} \left(\sum_{k=1}^Na_{ik}(y_{ij}-y_{kj})+a_{ij}(y_{ij}-x_j)\right),
\end{aligned}
\end{equation}
for $i,j\in\{1,2,\cdots,N\}$ and $\theta_{ij}=\theta\bar{\theta}_{ij}$, where $\theta$ is a positive parameter to be determined and $\bar{\theta}_{ij}$ is a fixed positive constant for each $i,j\in\{1,2,\cdots,N\}.$ Moreover, $\mathbf{y}_i=[y_{i1},y_{i2},\cdots,y_{iN}]^T$ stands for player $i$'s local estimate on $\mathbf{x}$ and $\nabla_if_i(\mathbf{y}_i)$ is defined as $\nabla_if_i(\mathbf{y}_i)=\frac{\partial f_i(\mathbf{x})}{\partial x_i}|_{\mathbf{x}=\mathbf{y}_i}.$ Furthermore, $a_{ij}$ is the element on the $i$th row and $j$th column of the adjacency matrix of $\mathcal{G}$.

Then, the concatenated-vector form of \eqref{bound_dis} is
\begin{equation}
\begin{aligned}
\dot{\mathbf{x}}=&-\left[\rho_{\bar{U}}\left(\nabla_if_i(\mathbf{y}_i)\right)\right]_{vec}\\
\dot{\mathbf{y}}=&-\theta \bar{\Theta}(\mathcal{L}\otimes I_{N\times N}+\mathcal{A})(\mathbf{y}-\mathbf{1}_N\otimes \mathbf{x}),
\end{aligned}
\end{equation}
where $\mathbf{y}=[y_{ij}]_{vec}, \bar{\Theta}=\text{diag}\{\bar{\theta}_{ij}\},$ $\mathcal{L}$ is the Laplacian matrix of $\mathcal{G}$, $\mathcal{A}=\text{diag}\{a_{ij}\}$ and $I_{N\times N}$ is an $N\times N$ dimensional identity matrix.


The following theorem establishes the stability result for the seeking strategy in \eqref{bound_dis}.

\begin{Theorem}\label{th1}
Suppose that Assumptions \ref{Assu_1}-\ref{Ass4} are satisfied, and the players update their actions according to \eqref{bound_dis}. Then, there exists a $\theta^*$ such that for each $\theta\in (\theta^*,\infty)$, the Nash equilibrium is globally asymptotically stable.
\end{Theorem}
\begin{Proof}
Define the Lyapunov candidate function as
\begin{equation}
V=\sum_{i=1}^N \int_0^{\nabla_if_i(\mathbf{x})} \rho_{\bar{U}}(t)dt+(\mathbf{y}-\mathbf{1}_N\otimes \mathbf{x})^T\mathcal{P} (\mathbf{y}-\mathbf{1}_N\otimes \mathbf{x})
\end{equation}
where $\mathcal{P}$ is a symmetric positive definite matrix that satisfies $\mathcal{P}\bar{\Theta} (\mathcal{L}\otimes I_{N\times N}+\mathcal{A})+(\mathcal{L}\otimes I_{N\times N}+\mathcal{A})\bar{\Theta}\mathcal{P}=\mathcal{Q}$
where $\mathcal{Q}$ is a symmetric positive definite matrix by Assumption \ref{a3} \cite{YETcyber18}. Then,
\begin{equation}
\begin{aligned}
\dot{V}=&-\left[\rho_{\bar{U}}\left(\nabla_if_i(\mathbf{x})\right)\right]_{vec}^TH(\mathbf{x})\left[\rho_{\bar{U}}\left(\nabla_if_i(\mathbf{y}_i)\right)\right]_{vec}\\
&+(\dot{\mathbf{y}}-\mathbf{1}_N\otimes \dot{\mathbf{x}})^T\mathcal{P}(\mathbf{y}-\mathbf{1}_N\otimes \mathbf{x})\\
&+(\mathbf{y}-\mathbf{1}_N\otimes \mathbf{x})^T\mathcal{P} (\dot{\mathbf{y}}-\mathbf{1}_N\otimes \dot{\mathbf{x}})\\
\leq &-\left[\rho_{\bar{U}}\left(\nabla_if_i(\mathbf{x})\right)\right]_{vec}^T H(\mathbf{x})\left[\rho_{\bar{U}}\left(\nabla_if_i(\mathbf{y}_i)\right)\right]_{vec}\\
&+2(\mathbf{y}-\mathbf{1}_N\otimes \mathbf{x})^T \mathcal{P}\left(\mathbf{1}_N\otimes \left[\rho_{\bar{U}}\left(\nabla_if_i(\mathbf{y}_i)\right)\right]_{vec}\right)\\
&-\lambda_{min}(\mathcal{Q})\theta||\mathbf{y}-\mathbf{1}_N\otimes \mathbf{x}||^2,
\end{aligned}
\end{equation}
in which $-\left[\rho_{\bar{U}}\left(\nabla_if_i(\mathbf{x})\right)\right]_{vec}^T H(\mathbf{x})\left[\rho_{\bar{U}}\left(\nabla_if_i(\mathbf{y}_i)\right)\right]_{vec}=-\left[\rho_{\bar{U}}\left(\nabla_if_i(\mathbf{x})\right)\right]_{vec}^TH(\mathbf{x})\left[\rho_{\bar{U}}\left(\nabla_if_i(\mathbf{x})\right)\right]_{vec}+\left[\rho_{\bar{U}}\left(\nabla_if_i(\mathbf{x})\right)\right]_{vec}^TH(\mathbf{x})
\left[\rho_{\bar{U}}\left(\nabla_if_i(\mathbf{x})\right)-\rho_{\bar{U}}\left(\nabla_if_i(\mathbf{y}_i)\right)\right]_{vec}.$
Furthermore, by Assumption \ref{Ass4}
$\left|\rho_{\bar{U}}\left(\nabla_if_i(\mathbf{x})\right)-\rho_{\bar{U}}\left(\nabla_if_i(\mathbf{y}_i)\right)\right|
\leq \left|\nabla_if_i(\mathbf{x})-\nabla_if_i(\mathbf{y}_i)\right|\leq \bar{l}_i||\mathbf{x}-\mathbf{y}_i||$, where $\bar{l}_i$ is the Lipschitz constant of $\nabla_if_i(\mathbf{x}).$
Noticing that the elements in $H(\mathbf{x})$ are bounded according to Assumption \ref{Ass4}, let $l_1=\sup_{\mathbf{x}\in \mathbb{R}^N} ||H(\mathbf{x})||\max\{\bar{l}_i\},$
$l_2=2||\mathcal{P}||\sqrt{N}\max\{\bar{l}_i\},$ and $l_3=2||\mathcal{P}||\sqrt{N},$ in which $\max\{\bar{l}_i\}$ is the maximum value of $\bar{l}_i$ for $i\in\{1,2,\cdots,N\}.$
Then, $\dot{V}\leq -m \left|\left|\left[\rho_{\bar{U}}\left(\nabla_if_i(\mathbf{x})\right)\right]_{vec}\right|\right|^2-\lambda_{min}(\mathcal{Q})\theta||\mathbf{y}-\mathbf{1}_N\otimes \mathbf{x}||^2+l_1 \left|\left|\left[\rho_{\bar{U}}\left(\nabla_if_i(\mathbf{x})\right)\right]_{vec}\right|\right| ||\mathbf{y}-\mathbf{1}_N\otimes \mathbf{x}||+l_2||\mathbf{y}-\mathbf{1}_N\otimes \mathbf{x}||^2+l_3||\mathbf{y}-\mathbf{1}_N\otimes \mathbf{x}||\left|\left|\left[\rho_{\bar{U}}\left(\nabla_if_i(\mathbf{x})\right)\right]_{vec}\right|\right|.$

Moreover, as $l_1 \left|\left|\left[\rho_{\bar{U}}\left(\nabla_if_i(\mathbf{x})\right)\right]_{vec}\right|\right| ||\mathbf{y}-\mathbf{1}_N\otimes \mathbf{x}||\leq \frac{l_1}{2\epsilon_1}\left|\left| \left[\rho_{\bar{U}}\left(\nabla_if_i(\mathbf{x})\right)\right]_{vec}\right|\right|^2+\frac{l_1\epsilon_1}{2}||\mathbf{y}-\mathbf{1}_N\otimes \mathbf{x}||^2$ and
$l_3||\mathbf{y}-\mathbf{1}_N\otimes \mathbf{x}||\left|\left|\left[\rho_{\bar{U}}\left(\nabla_if_i(\mathbf{x})\right)\right]_{vec}\right|\right|\leq \frac{l_3}{2\epsilon_2}\left|\left| \left[\rho_{\bar{U}}\left(\nabla_if_i(\mathbf{x})\right)\right]_{vec}\right|\right|^2+\frac{l_3\epsilon_2}{2}||\mathbf{y}-\mathbf{1}_N\otimes \mathbf{x}||^2,$
in which $\epsilon_1,\epsilon_2$ are positive constants that can be arbitrarily chosen,
\begin{equation}
\begin{aligned}
&\dot{V}\leq -\left(m-\frac{l_1}{2\epsilon_1}-\frac{l_3}{2\epsilon_2}\right)\left|\left| \left[\rho_{\bar{U}}\left(\nabla_if_i(\mathbf{x})\right)\right]_{vec}\right|\right|^2\\
&-\left(\lambda_{min}(\mathcal{Q})\theta -l_2-\frac{l_1\epsilon_1}{2}-\frac{l_3\epsilon_2}{2}\right)||\mathbf{y}-\mathbf{1}_N\otimes \mathbf{x}||^2.
\end{aligned}
\end{equation}

Choose $\epsilon_1,\epsilon_2$ such that $m-\frac{l_1}{2\epsilon_1}-\frac{l_3}{2\epsilon_2}>0$. Then, for fixed $\epsilon_1,\epsilon_2$, let
\begin{equation}\label{lo_b}
\theta^* =\frac{2l_2+l_1\epsilon_1+l_3\epsilon_2}{2\lambda_{min}(\mathcal{Q})},
\end{equation}
and $\theta>\theta^*.$ Subsequently, let $l_4=\min\{m-\frac{l_1}{2\epsilon_1}-\frac{l_3}{2\epsilon_2},\lambda_{min}(\mathcal{Q})\theta -l_2-\frac{l_1\epsilon_1}{2}-\frac{l_3\epsilon_2}{2}\}$, then, $\dot{V}\leq -l_4 ||\mathbf{\chi}||^2,$
where $\mathbf{\chi}=\left[\left[\rho_{\bar{U}}\left(\nabla_if_i(\mathbf{x})\right)\right]_{vec}^T,(\mathbf{y}-\mathbf{1}_N\otimes \mathbf{x})^T\right]^T.$
Moreover, following the proof of Theorem \ref{Lem_1}, it can be shown that $V$ is positive definite and radially unbounded with respect to $\mathbf{\chi}$. Hence, $||\mathbf{\chi}||\rightarrow 0$ as $t\rightarrow +\infty.$ Recalling that by Assumption \ref{Assu_2}, $\bar{\mathcal{P}}(\mathbf{x})=\mathbf{0}_N$ if and only if $\mathbf{x}=\mathbf{x}^*$, we see that   $\mathbf{y}\rightarrow \mathbf{1}_N\otimes \mathbf{x}\rightarrow \mathbf{1}_N\otimes \mathbf{x}^*$ as $t\rightarrow +\infty.$  To this end, we arrive at the conclusion.
\end{Proof}
\begin{Remark}
The seeking strategy in \eqref{bound_dis} is adapted from the seeking strategy in \cite{YETAC17} in which the saturation function is included to ensure that $|u_i|\leq \bar{U}.$
\end{Remark}

\subsection{Second-order integrator-type systems}
In this section, we consider Nash equilibrium seeking for games in second-order integrator-type systems in which player $i$'s action is governed by
\begin{equation}\label{second}
\dot{x}_i=\nu_i,\dot{\nu}_i=u_i,
\end{equation}
for $i\in\{1,2,\cdots,N\}.$
More specifically, in Section \ref{ce_2_1}, a centralized algorithm will be proposed without considering the boundedness of the control inputs. Moreover, the problem is reconsidered under distributed networks in Section \ref{ce_2_2}. Lastly, the boundedness of the control inputs will be addressed in Section \ref{ce_2_3}.
\subsubsection{Centralized Nash equilibrium seeking without considering the boundedness of the control inputs}\label{ce_2_1}
Let the Nash equilibrium seeking strategy be
\begin{equation}\label{st_1}
\dot{\mathbf{x}}=\mathbf{\nu},\dot{\mathbf{\nu}}=-\alpha \small{\left[\nabla_if_i(\mathbf{x})\right]}_{vec}-\beta \mathbf{\nu}-H(\mathbf{x})\mathbf{\nu},
\end{equation}
where $\mathbf{\nu}=[\nu_i]_{vec}$ and $\alpha,\beta$ are positive control gains to be determined. Then, the following result can be obtained.

\begin{Theorem}\label{TH2}
Suppose that  Assumptions \ref{Assu_1} and \ref{Assu_2} are satisfied and the players update their actions according to \eqref{st_1}. Then, there exists a positive constant $\alpha^*$ such that for each $\alpha\in (0,\alpha^*),$ there exists a positive constant $\beta^*(\alpha)$ such that for each $\beta \in(0,\beta^*),$ the Nash equilibrium is globally asymptotically stable under \eqref{st_1}.
\end{Theorem}
\begin{Proof}
Define the Lyapunov candidate function as
\begin{equation}
V=\mathbf{\nu}^T\mathbf{\nu}+\frac{1}{2}\small{\left[\nabla_if_i(\mathbf{x})\right]_{vec}^T\left[\nabla_if_i(\mathbf{x})\right]_{vec}}+\mathbf{\nu}^T\small{\left[\nabla_if_i(\mathbf{x})\right]_{vec}}.
\end{equation}
Then, $V=\frac{1}{6}\left|\left|\left[\nabla_if_i(\mathbf{x})\right]_{vec}\right|\right|^2+\frac{1}{4}||\mathbf{\nu}||^2+\left|\left|\frac{1}{\sqrt{3}}\left[\nabla_if_i(\mathbf{x})\right]_{vec}+\frac{\sqrt{3}}{2}\mathbf{\nu}\right|\right|^2,$
and it can be easily concluded that the Lyapunov candidate function is positive definite and radially unbounded with respect to $\left[\mathbf{\nu}^T,\left[\nabla_if_i(\mathbf{x})\right]_{vec}^T\right]^T$. Moreover, by Assumption \ref{Assu_2},
\begin{equation}
\begin{aligned}
\dot{V}=&2\mathbf{\nu}^T\left(-\beta\mathbf{\nu}-\alpha\left[\nabla_if_i(\mathbf{x})\right]_{vec}-H(\mathbf{x})\mathbf{\nu}\right)\\
&+\left[\nabla_if_i(\mathbf{x})\right]_{vec}^TH(\mathbf{x})\mathbf{\nu}+\mathbf{\nu}^TH(\mathbf{x})\mathbf{\nu}\\
&+\left(-\beta \mathbf{\nu}-\alpha\left[\nabla_if_i(\mathbf{x})\right]_{vec}-H(\mathbf{x})\mathbf{\nu}\right)^T\left[\nabla_if_i(\mathbf{x})\right]_{vec}\\
\leq & -(2\beta+m)||\mathbf{\nu}||^2-\alpha \left|\left|\left[\nabla_if_i(\mathbf{x})\right]_{vec}\right|\right|^2\\
&+(2\alpha+\beta)||\mathbf{\nu}||\left|\left|\left[\nabla_if_i(\mathbf{x})\right]_{vec}\right|\right|\\
\leq & -\left(2\beta+m-(2\alpha+\beta)/(2\epsilon_1)\right)||\mathbf{\nu}||^2\\
&-\left(\alpha-(\epsilon_1(2\alpha+\beta))/2\right)\left|\left|\left[\nabla_if_i(\mathbf{x})\right]_{vec}\right|\right|^2,
\end{aligned}
\end{equation}
where $\epsilon_1$ is a positive constant that can be arbitrarily chosen. Let $\frac{2\alpha+\beta}{2(2\beta+m)}<\epsilon_1<\frac{2\alpha}{2\alpha+\beta},2\alpha-2\sqrt{\alpha m}<\beta<2\alpha+2\sqrt{\alpha m}.$
Then, $\dot{V}$ is negative definite. Hence, the conclusion can be drawn with $\alpha^*=m$ and $\beta^*=2\alpha+2\sqrt{\alpha m}.$
\end{Proof}

The seeking strategy in \eqref{st_1} achieves the Nash equilibrium seeking in a centralized fashion. However, as it is challenging for the players to simultaneously estimate $H(\mathbf{x})$ and $\mathbf{x}$ in a distributed fashion, the distributed implementation of \eqref{st_1} is difficult to be achieved. In the following, we consider Nash equilibrium seeking in distributed networks from another perspective.

\subsubsection{Distributed Nash equilibrium seeking without considering the boundedness of control inputs}\label{ce_2_2}
Suppose that in the considered game, each player $i,i\in\{1,2,\cdots,N\}$ updates their own action according to
\begin{equation}\label{method_3f}
\begin{aligned}
&\dot{x}_i=\nu_i,\dot{\nu}_i=-(x_i-z_i)-(\nu_i-\dot{z}_i),\dot{z}_i=-\bar{K}_i\nabla_if_i(\mathbf{y}_i)\\
&\dot{y}_{ij}=-\theta_{ij}(\sum_{k=1}^N a_{ik}(y_{ij}-y_{kj})+a_{ij}(y_{ij}-z_j)),
\end{aligned}
\end{equation}
where $j\in\{1,2,\cdots,N\}$ and $z_i,y_{ij}$ are auxiliary variables. Moreover, $\bar{K}_i=\theta_1K_i$, $\theta_{ij}=\theta\theta_1\bar{\theta}_{ij}$ in which $\theta,\theta_1$ are positive parameters to be determined and $K_i,\bar{\theta}_{ij}$ are fixed positive constants.

The concatenated vector form of \eqref{method_3f} is
\begin{equation}\label{method_3}
\begin{aligned}
\dot{\mathbf{x}}&=\mathbf{\nu},\dot{\mathbf{\nu}}=-(\mathbf{x}-\mathbf{z})-(\mathbf{\nu}-\dot{\mathbf{z}})\\
\dot{\mathbf{z}}&=-\bar{K} \small{\left[\nabla_if_i(\mathbf{y}_i)\right]_{vec}}\\
\dot{\mathbf{y}}&=-\Theta(\mathcal{L}\otimes I_{N\times N}+\mathcal{A})(\mathbf{y}-\mathbf{1}_N\otimes \mathbf{z}),
\end{aligned}
\end{equation}
where $\bar{K}=\text{diag}\{\bar{K}_i\}, \Theta=\text{diag}\{\theta_{ij}\}$ and $\mathbf{z}=[z_i]_{vec}$.

The following theorem establishes the stability of the equilibrium in \eqref{method_3}.

\begin{Theorem}\label{th3}
Suppose that Assumptions 1-4 are satisfied and the players update their actions according to \eqref{method_3}. Then, there exists a positive constant $\theta^*$ such that for each $\theta\in(\theta^*,\infty)$, there exists a positive constant $\theta_1^*(\theta)$ such that for each $\theta_1\in (0,\theta_1^*),$ the Nash equilibrium is globally asymptotically stable.
\end{Theorem}
\begin{Proof}
Consider
\begin{equation}
\begin{aligned}
V(\mathbf{\eta})=&\frac{1}{2}(\mathbf{z}-\mathbf{x}^*)^TK^{-1}(\mathbf{z}-\mathbf{x}^*)\\
&+(\mathbf{y}-\mathbf{1}_N\otimes \mathbf{z})^T\mathcal{P}(\mathbf{y}-\mathbf{1}_N\otimes \mathbf{z})\\
&+\frac{1}{2}(\mathbf{x}-\mathbf{z})^T(\mathbf{x}-\mathbf{z})+\frac{1}{2}(\mathbf{\nu}-\dot{\mathbf{z}})^T(\mathbf{\nu}-\dot{\mathbf{z}}),
\end{aligned}
\end{equation}
where $\mathcal{P}$ is defined in the proof of Theorem \ref{th1}, $\mathbf{\eta}=[(\mathbf{z}-\mathbf{x}^*)^T,(\mathbf{y}-\mathbf{1}_N\otimes \mathbf{z})^T,(\mathbf{x}-\mathbf{z})^T,(\mathbf{\nu}-\dot{\mathbf{z}})^T]^T$ and $K=\text{diag}\{K_i\}$ as the Lyapunov candidate function.
Then,
\begin{equation}
\begin{aligned}
\dot{V}\leq &-\theta_1 (\mathbf{z}-\mathbf{x}^*)^T\left[\nabla_if_i(\mathbf{z})\right]_{vec}\\
&-\lambda_{min}(\mathcal{Q})\theta \theta_1 ||\mathbf{y}-\mathbf{1}_N\otimes \mathbf{z}||^2\\
&-||\mathbf{\nu}-\dot{\mathbf{z}}||^2+\theta_1(\mathbf{z}-\mathbf{x}^*)^T\left[\nabla_if_i(\mathbf{z})-\nabla_if_i(\mathbf{y}_i)\right]_{vec}\\
&-2(\mathbf{y}-\mathbf{1}_N\otimes \mathbf{z})^T\mathcal{P}\mathbf{1}_N\otimes \dot{\mathbf{z}}-(\mathbf{\nu}-\dot{\mathbf{z}})^T \ddot{\mathbf{z}}.
\end{aligned}
\end{equation}
By Assumption \ref{Assu_2}, $-(\mathbf{z}-\mathbf{x}^*)^T\left[\nabla_if_i(\mathbf{z})\right]_{vec}\leq -m||\mathbf{z}-\mathbf{x}^*||^2.$
Moreover, by Assumption \ref{Assu_1}, there exists positive constant $\bar{l}_{i}$ such that $\left|\left|\nabla_if_i(\mathbf{z})-\nabla_if_i(\mathbf{y}_i)\right|\right|\leq \bar{l}_{i}||\mathbf{y}-\mathbf{1}_N\otimes \mathbf{z}||,$
and $\left|\left|\nabla_if_i(\mathbf{y}_i)\right|\right|=\left|\left|\nabla_if_i(\mathbf{y}_i)-\nabla_if_i(\mathbf{z})+\nabla_if_i(\mathbf{z})-\nabla_if_i(\mathbf{x}^*)\right|\right|
\leq \bar{l}_{i}||\mathbf{y}_i-\mathbf{z}||+\bar{l}_{i}||\mathbf{z}-\mathbf{x}^*||.$

In addition, $\ddot{\mathbf{z}}=\theta \theta_1^2 K\bar{H}(\mathbf{y})\bar{\Theta}(\mathcal{L}\otimes I_{N\times N}+\mathcal{A})(\mathbf{y}-\mathbf{1}_N\otimes\mathbf{z}),$
where $\bar{H}(\mathbf{y})=\left[
                             \begin{array}{cccc}
                               \bar{h}_{11} & \bar{h}_{12} & \cdots & \bar{h}_{1N} \\
                               \bar{h}_{21} & \bar{h}_{22} & \cdots & \bar{h}_{2N} \\
                               \vdots &  & \ddots &  \\
                               \bar{h}_{N1} & \bar{h}_{N2} & \cdots & \bar{h}_{NN} \\
                             \end{array}
                           \right]
$ and $\bar{h}_{ij}\in \mathbb{R}^{1\times N}$. Moreover, $\bar{h}_{ij}=\mathbf{0}^T_{N}$ for $i\neq j$ and $\bar{h}_{ii}=\left[\frac{\partial^2 f_i}{\partial x_i\partial x_1}(\mathbf{y}_i),\frac{\partial^2 f_i}{\partial x_i\partial x_2}(\mathbf{y}_i),\cdots,\frac{\partial^2 f_i}{\partial x_i\partial x_N}(\mathbf{y}_i)\right],$ where $\frac{\partial^2 f_i}{\partial x_i\partial x_j}(\mathbf{y}_i)=\frac{\partial^2 f_i(\mathbf{x})}{\partial x_i\partial x_j}\left.\right|_{\mathbf{x}=\mathbf{y}_i}.$ Noticing that $\bar{H}(\mathbf{y})$ is bounded according to Assumption 4, let $l_1=\max\{\bar{l}_i\}+2||\mathcal{P}||N\max\{K_i\bar{l}_i\},l_2=2||\mathcal{P}||\sqrt{N}\max\{K_i\bar{l}_i\}$ and $l_3=||K||\sup_{\mathbf{y}}||\bar{H}(\mathbf{y})||||\bar{\Theta}(\mathcal{L}\otimes I_{N\times N}+\mathcal{A})||$. Then,
\begin{equation}
\begin{aligned}
\dot{V}\leq &-\theta_1m||\mathbf{z}-\mathbf{x}^*||^2-\lambda_{min}(\mathcal{Q})\theta\theta_1||\mathbf{y}-\mathbf{1}_N\otimes \mathbf{z}||^2\\
&-||\mathbf{\nu}-\dot{\mathbf{z}}||^2+\theta_1l_1||\mathbf{z}-\mathbf{x}^*||||\mathbf{y}-\mathbf{1}_N\otimes \mathbf{z}||\\
&+\theta_1 l_2||\mathbf{y}-\mathbf{1}_N\otimes \mathbf{z}||^2+\theta\theta_1^2l_3||\mathbf{\nu}-\dot{\mathbf{z}}||||\mathbf{y}-\mathbf{1}_N\otimes \mathbf{z}||.
\end{aligned}
\end{equation}
Define $A_1=\left[
         \begin{array}{cc}
           m & -\frac{l_1}{2} \\
           -\frac{l_1}{2} & \lambda_{min}(\mathcal{Q})\theta-l_2 \\
         \end{array}
       \right],$ and choose $\theta>\theta^*$ where
       \begin{equation}\label{b1}
       \theta^*=\frac{l_1^2}{4m\lambda_{min}(\mathcal{Q})}+\frac{l_2}{\lambda_{min}(\mathcal{Q})},
       \end{equation}
        then, $\dot{V}\leq -\theta_1 \lambda_{min}(A_1)||E_1||^2-||\mathbf{\nu}-\dot{\mathbf{z}}||^2+\theta\theta_1^2l_3||\mathbf{\nu}-\dot{\mathbf{z}}||||\mathbf{y}-\mathbf{1}_N\otimes \mathbf{z}||,$
where $\lambda_{min}(A_1)>0$ and $E_1=[(\mathbf{z}-\mathbf{x}^*)^T,(\mathbf{y}-\mathbf{1}_N\otimes\mathbf{z})^T]^T.$
Moreover, define $A_2=\left[
         \begin{array}{cc}
           \theta_1 \lambda_{min}(A_1) & -\frac{\theta \theta_1^2l_3}{2} \\
           -\frac{\theta \theta_1^2l_3}{2} & 1 \\
         \end{array}
       \right].$ Then, $\lambda_{min}(A_2)>0$ given that $\theta_1<\theta_1^*$, where
       \begin{equation}\label{b2}
       \theta_1^*=\small{\left(\frac{4\lambda_{min}(A_1)}{\theta^2 l_3^2}\right)^{\frac{1}{3}}}.
       \end{equation}

If this is the case,
\begin{equation}
\dot{V}\leq -\lambda_{min}(A_2)||E||^2,
\end{equation}
where $E=[(\mathbf{z}-\mathbf{x}^*)^T,(\mathbf{y}-\mathbf{1}_N\otimes \mathbf{z})^T,(\mathbf{\nu}-\dot{\mathbf{z}})^T]^T.$
Hence, $\mathbf{z}=\mathbf{x}^*,\mathbf{y}=\mathbf{1}_N\otimes \mathbf{z},\mathbf{\nu}=\dot{\mathbf{z}}$ at $\dot{V}=0,$ which indicates that $\dot{\mathbf{x}}=\mathbf{\nu},\dot{\mathbf{\nu}}=-(\mathbf{x}-\mathbf{z}),\dot{\mathbf{z}}=\mathbf{0}_N,\dot{\mathbf{y}}=\mathbf{0}_{N^2}.$
Recalling that $\mathbf{\nu}=\dot{\mathbf{z}}$ at $\dot{V}=0$, we have $\mathbf{\nu}=\mathbf{0}$. Hence, $\dot{\mathbf{x}}=\mathbf{0}_N$ and $\mathbf{x}=C_1,\mathbf{z}=C_2$ at $\dot{V}=0,$ where $C_1,C_2$ are constant vectors. Therefore, $\dot{\mathbf{\nu}}=-C_1+C_2,$
at $\dot{V}=0.$ Recalling that $\mathbf{\nu}=0$, we can get that $C_1=C_2$, i.e., $\mathbf{x}=\mathbf{z}.$ Hence, the conclusion can be derived by utilizing the LaSalle's invariance principle.
\end{Proof}

The strategy in \eqref{method_3} addressed the Nash equilibrium seeking problem for games in second-order integrator-type systems without considering the boundedness of the controls. In the upcoming section, the seeking strategy  in \eqref{method_3} will be adapted for systems where the controls are bounded.

\subsubsection{Distributed Nash equilibrium with bounded control inputs}\label{ce_2_3}
Let the Nash equilibrium seeking strategy be
\begin{equation}\label{method_4}
\begin{aligned}
\dot{\mathbf{x}}&=\mathbf{\nu},\dot{\mathbf{\nu}}=-\rho_{\bar{U}}((\mathbf{x}-\mathbf{z})+(\mathbf{\nu}-\dot{\mathbf{z}}))\\
\dot{\mathbf{z}}&=-\bar{K} \small{\left[\nabla_if_i(\mathbf{y}_i)\right]_{vec}}\\
\dot{\mathbf{y}}&=-\Theta(\mathcal{L}\otimes I_{N\times N}+\mathcal{A})(\mathbf{y}-\mathbf{1}_N\otimes \mathbf{z}).
\end{aligned}
\end{equation}
Then, the following result can be derived.
\begin{Theorem}\label{th4}
Suppose that Assumptions \ref{Assu_1}-\ref{Assu_2} are satisfied. Then, for any positive constant $\Delta$, there exists a positive constant $\theta^*$ such that for each $\theta\in(\theta^*,\infty)$, there exists a positive constant $\theta_1^*(\Delta,\theta)$ such that for each $\theta_1\in(0,\theta_1^*)$, $\mathbf{x}$ generated by \eqref{method_4} converges asymptotically to $\mathbf{x}^*$ given that $||(\mathbf{\nu}(0)-\dot{\mathbf{z}}(0))^T,(\mathbf{x}(0)-\mathbf{z}(0))^T,(\mathbf{y}(0)-\mathbf{1}_N\otimes \mathbf{z}(0))^T,\left(\mathbf{z}(0)-\mathbf{x}^*\right)^T||\leq \Delta.$
\end{Theorem}
\begin{Proof}
Define the Lyapunov candidate function as
\begin{equation}
\begin{aligned}
&V(\mathbf{\eta})=\frac{1}{2}(\mathbf{z}-\mathbf{x}^*)^TK^{-1}(\mathbf{z}-\mathbf{x}^*)+(\mathbf{\nu}-\dot{\mathbf{z}})^T(\mathbf{\nu}-\dot{\mathbf{z}})\\
&+(\mathbf{y}-\mathbf{1}_N\otimes \mathbf{z})^T\mathcal{P}(\mathbf{y}-\mathbf{1}_N\otimes \mathbf{z})+\\
&\small{\sum_{i=1}^N \int_{0}^{x_i-z_i}\rho_{\bar{U}}(t)dt}+\small{\sum_{i=1}^N\int_{0}^{x_i-z_i+v_i-\dot{z}_i}\rho_{\bar{U}}(t)dt},
\end{aligned}
\end{equation}
where $\mathcal{P}$ is defined in the proof of Theorem 1 and $\mathbf{\eta}=[(\mathbf{z}-\mathbf{x}^*)^T,(\mathbf{y}-\mathbf{1}_N\otimes \mathbf{z})^T,(\mathbf{x}-\mathbf{z})^T,(\mathbf{\nu}-\dot{\mathbf{z}})^T,(\mathbf{x}-\mathbf{z}+\mathbf{\nu}-\dot{\mathbf{z}})^T]^T$. Then, it can be easily derived that the Lyapunov candidate function is positive definite and radially unbounded.
Moreover, following the analysis in the proof of Theorem \ref{th3}, it can be derived that
\begin{equation}
\begin{aligned}
\dot{V}\leq &-\theta_1m||\mathbf{z}-\mathbf{x}^*||^2-\lambda_{min}(\mathcal{Q})\theta \theta_1||\mathbf{y}-\mathbf{1}_N\otimes \mathbf{z}||^2\\
&+\theta_1l_1||\mathbf{z}-\mathbf{x}^*||||\mathbf{y}-\mathbf{1}_N\otimes\mathbf{z}||+\theta_1l_2||\mathbf{y}-\mathbf{1}_N\otimes\mathbf{z}||^2\\
&+\rho_{\bar{U}}(\mathbf{x}-\mathbf{z})^T(\mathbf{\nu}-\dot{\mathbf{z}})-2(\mathbf{\nu}-\dot{\mathbf{z}})^T\rho_{\bar{U}}(\mathbf{x}-\mathbf{z}+\mathbf{\nu}-\dot{\mathbf{z}})\\
&-2(\mathbf{\nu}-\dot{\mathbf{z}})^T\ddot{\mathbf{z}}+\rho_{\bar{U}}(\mathbf{x}-\mathbf{z}+\mathbf{\nu}-\mathbf{\dot{z}})^T(\mathbf{\nu}-\dot{\mathbf{z}})\\
&-\rho_{\bar{U}}(\mathbf{x}-\mathbf{z}+\mathbf{\nu}-\dot{\mathbf{z}})^T\rho_{\bar{U}}(\mathbf{x}-\mathbf{z}+\mathbf{\nu}-\dot{\mathbf{z}})\\
&-\rho_{\bar{U}}(\mathbf{x}-\mathbf{z}+\mathbf{\nu}-\dot{\mathbf{z}})^T\ddot{\mathbf{z}},
\end{aligned}
\end{equation}
where $l_1=\max\{\bar{l}_i\}+2||\mathcal{P}||N\max\{K_i\bar{l}_i\}$ and $l_2=2||\mathcal{P}||\sqrt{N}\max\{K_i\bar{l}_i\}$.
Since $-(\mathbf{\nu}-\dot{\mathbf{z}})^T(\rho_{\bar{U}}(\mathbf{x}-\mathbf{z}+\mathbf{\nu}-\dot{\mathbf{z}})-\rho_{\bar{U}}(\mathbf{x}-\mathbf{z}))\leq 0,$
we have
$-(\mathbf{\nu}-\dot{\mathbf{z}})^T(\rho_{\bar{U}}(\mathbf{x}-\mathbf{z}+\mathbf{\nu}-\dot{\mathbf{z}})-\rho_{\bar{U}}(\mathbf{x}-\mathbf{z}))-\rho_{\bar{U}}(\mathbf{x}-\mathbf{z}+\mathbf{\nu}-\dot{\mathbf{z}})^T\rho_{\bar{U}}(\mathbf{x}-\mathbf{z}+\mathbf{\nu}-\dot{\mathbf{z}})=0$ if and only if $(\mathbf{\nu}-\dot{\mathbf{z}})^T(\rho_{\bar{U}}(\mathbf{x}-\mathbf{z}+\mathbf{\nu}-\dot{\mathbf{z}})-\rho_{\bar{U}}(\mathbf{x}-\mathbf{z}))=0$ and $\rho_{\bar{U}}(\mathbf{x}-\mathbf{z}+\mathbf{\nu}-\dot{\mathbf{z}})=\mathbf{0}.$ Moreover, from  $\rho_{\bar{U}}(\mathbf{x}-\mathbf{z}+\mathbf{\nu}-\dot{\mathbf{z}})=\mathbf{0},$ we have $\mathbf{x}-\mathbf{z}+\mathbf{\nu}-\dot{\mathbf{z}}=\mathbf{0},$ by which $(\mathbf{\nu}-\dot{\mathbf{z}})^T\rho_{\bar{U}}(\mathbf{x}-\mathbf{z})=0.$ Therefore, $(\mathbf{\nu}-\dot{\mathbf{z}})^T\rho_{\bar{U}}(\mathbf{\nu}-\dot{\mathbf{z}})=0,$ from which we can get that $\mathbf{\nu}-\dot{\mathbf{z}}=\mathbf{0}$ and $\mathbf{x}-\mathbf{z}=\mathbf{0}.$ Hence, $-(\mathbf{\nu}-\dot{\mathbf{z}})^T(\rho_{\bar{U}}(\mathbf{x}-\mathbf{z}+\mathbf{\nu}-\dot{\mathbf{z}})-\rho_{\bar{U}}(\mathbf{x}-\mathbf{z}))-\rho_{\bar{U}}(\mathbf{x}-\mathbf{z}+\mathbf{\nu}-\dot{\mathbf{z}})^T\rho_{\bar{U}}(\mathbf{x}-\mathbf{z}+\mathbf{\nu}-\dot{\mathbf{z}})$ is negative definite with respect to $[(\mathbf{x}-\mathbf{z})^T,(\mathbf{\nu}-\dot{\mathbf{z}})^T,(\mathbf{x}-\mathbf{z}+\mathbf{\nu}-\dot{\mathbf{z}})^T]^T.$

By further following the proof of Theorem \ref{th3}, we can conclude that by choosing $\theta>\theta^*$, where
\begin{equation}\label{b3}
\theta^*=\frac{l_1^2+4ml_2}{4m\lambda_{min}(\mathcal{Q})},
\end{equation}
we have $\dot{V}\leq -\theta_1\lambda_{min}(A_1)||E_1||^2-W_2(E_2)+l_3\theta \theta_1^2||\rho_{\bar{U}}(\mathbf{x}-\mathbf{z}+\mathbf{\nu}-\dot{\mathbf{z}})||||\mathbf{y}-\mathbf{1}_N\otimes \mathbf{z}||+2l_3 \theta \theta_1^2||\mathbf{\nu}-\dot{\mathbf{z}}||||\mathbf{y}-\mathbf{1}_N\otimes \mathbf{z}||,$
where $A_1=\left[
         \begin{array}{cc}
           m & -\frac{l_1}{2} \\
           -\frac{l_1}{2} & \lambda_{min}(\mathcal{Q})\theta-l_2 \\
         \end{array}
       \right]$, $E_1=[(\mathbf{z}-\mathbf{x}^*)^T,(\mathbf{y}-\mathbf{1}_N\otimes \mathbf{z})^T]^T,$ $W_2(E_2)=(\mathbf{\nu}-\dot{\mathbf{z}})^T(\rho_{\bar{U}}(\mathbf{x}-\mathbf{z}+\mathbf{\nu}-\dot{\mathbf{z}})-\rho_{\bar{U}}(\mathbf{x}-\mathbf{z}))+\rho_{\bar{U}}(\mathbf{x}-\mathbf{z}+\mathbf{\nu}-\dot{\mathbf{z}})^T\rho_{\bar{U}}(\mathbf{x}-\mathbf{z}+\mathbf{\nu}-\dot{\mathbf{z}}),$ $E_2=[(\mathbf{x}-\mathbf{z})^T,(\mathbf{\nu}-\dot{\mathbf{z}})^T,(\mathbf{x}-\mathbf{z}+\mathbf{\nu}-\dot{\mathbf{z}})^T]^T,$ and
       $l_3=||K||\sup_{\mathbf{y}}||\bar{H}(\mathbf{y})||||\bar{\Theta}(\mathcal{L}\otimes I_{N\times N}+\mathcal{A})||.$

To facilitate the subsequent analysis, define $W(\mathbf{\eta})=\lambda_{min}(A_1)||E_1||^2+W_2(E_2).$
Then, it is clear that $W(\mathbf{\eta})$ is positive definite and there exists a class $\mathcal{K}$ function $\gamma$ such that $\gamma(||\mathbf{\eta}||)\leq W(\mathbf{\eta}).$ Hence, if we choose $\theta_1\leq 1$, one can obtain that $-\theta_1\lambda(A_1)||E_1||^2-W_2(E_2)\leq -\theta_1W(\mathbf{\eta})\leq -\theta_1\gamma(||\mathbf{\eta}||).$ Similarly, if $\theta_1>1,$ one has $-\theta_1\lambda(A_1)||E_1||^2-W_2(E_2)\leq -\gamma(||\mathbf{\eta}||).$ Therefore,
$\dot{V}\leq -\min\{\theta_1,1\}\gamma(||\mathbf{\eta}||)+2l_3\theta \theta_1^2 ||\mathbf{\nu}-\dot{\mathbf{z}}||||\mathbf{y}-\mathbf{1}_N\otimes\mathbf{z}||+l_3\theta \theta_1^2||\rho_{\bar{U}}(\mathbf{x}-\mathbf{z}+\mathbf{\nu}-\dot{\mathbf{z}})||||\mathbf{y}-\mathbf{1}_N\otimes \mathbf{z}||.$


Therefore, for $\mathbf{\eta}$ that belongs to any compact set $D_{\omega}$ that contains the origin, $\dot{V}\leq -\min\{\theta_1,1\}\gamma(||\mathbf{\eta}||)+\theta \theta_1^2 l_4,$ where $l_4=\sup_{\mathbf{\eta}\in D_{\omega}} (2l_3||\mathbf{\nu}-\dot{\mathbf{z}}||||\mathbf{y}-\mathbf{1}_N\otimes \mathbf{z}||+l_3||\rho_{\bar{U}}(\mathbf{x}-\mathbf{z}+\mathbf{\nu}-\dot{\mathbf{z}})||||\mathbf{y}-\mathbf{1}_N\otimes \mathbf{z}||).$
Hence, $\dot{V}\leq -\frac{\min\{\theta_1,1\}}{2}\gamma(||\mathbf{\eta}||),\forall ||\mathbf{\eta}||\geq \gamma^{-1}(\frac{2\theta\theta_1^2l_4}{\min\{\theta_1,1\}}).$

Recalling that $V$ is positive definite, there exist $\gamma_1,\gamma_2\in \mathcal{K}$ such that $\gamma_1(||\mathbf{\eta}||)\leq V(\mathbf{\eta})\leq \gamma_2(||\mathbf{\eta}||).$
Take a positive constant $r$ such that $B_r\subset D_{\omega},$ where $B_r$ denotes an origin-centered ball with radius $r$. Moreover, choose $\theta_1$ to be sufficiently small such that $\gamma^{-1}(\frac{2\theta\theta_1^2l_4}{\min\{\theta_1,1\}})<\gamma_2^{-1}(\gamma_1(r))$. Then for any initial condition that satisfies $||\mathbf{\eta}(0)||\leq \gamma_2^{-1}(\gamma_1(r)),$ (i.e., $\Delta=\gamma_2^{-1}(\gamma_1(r))$), there exists a positive constant $T_1$ such that $||\mathbf{\eta}(t)||\leq \gamma_1^{-1}(\gamma_2(\gamma^{-1}(\frac{2\theta\theta_1^2l_4}{\min\{\theta_1,1\}})))$ for all $t\geq T_1$. Choosing $\theta_1$ to be sufficiently small such that $\gamma_1^{-1}(\gamma_2(\gamma^{-1}(\frac{2\theta\theta_1^2l_4}{\min\{\theta_1,1\}})))<\bar{U}$, then, the trajectory of the system in \eqref{method_4} is the same as the trajectory of the system in \eqref{method_3} for $t\geq T_1$ (with the same initial condition at $t=T_1$). Hence, further following the result in Theorem 4, the conclusion can be derived.
\end{Proof}

\begin{Remark}
Theorem \ref{th4} demonstrates a semi-global convergence result. That is, for any bounded initial conditions, the proposed method can drive the players' actions to the Nash equilibrium of the game by suitably tuning the control gains (possibly depend on the initial values of the variables). Different from local convergence results that require the initial errors to be sufficiently small, the semi-global results only require the initial values  to be bounded and the bounds can be arbitrarily large.
\end{Remark}
\begin{Remark}
Theorem \ref{th1} indicates that $\theta$ should be sufficiently large to ensure the convergence of \eqref{bound_dis}. The lower bound $\theta^*$ is qualified in \eqref{lo_b}. Similarly, Theorems
\ref{th3}-\ref{th4} illustrate that $\theta$ should be chosen to be sufficiently large while $\theta_1$ should be chosen to be sufficiently small to ensure the convergence of \eqref{method_3} and \eqref{method_4}. In the proof of the theorems, the lower bound of $\theta$ and upper bound $\theta_1$ are provided in \eqref{b1}-\eqref{b2} and \eqref{b3}, except that $\theta_1^*$ in Theorem \ref{th4} depends on the initial errors and is hard to be explicitly quantified without knowing the initial errors. From the quantifications of $\theta^*$ and $\theta_1^*$, it is clear that they depend on the Lipschitz constants and strong monotonicity constant of the pseudo-gradient vector, the communication topology as well as the number of players in the game. Moreover, though $\theta_1^*$ in  Theorem \ref{th4} depends on the unknown initial errors, the result is still meaningful as it suggests that for any bounded initial errors, we can directly choose $\theta_1$ to be sufficiently small to ensure the convergence of the proposed method. Interested readers are referred to the proofs of the corresponding theorems for more details.
\end{Remark}

\begin{Remark}
The theoretical results presented in the paper are established for $x_i\in \mathbb{R}.$ However, they can be easily extended to the case in which $x_i\in \mathbb{R}^p$ and $p\geq 2$ is a positive integer. Moreover, we suppose that the control inputs satisfy $|u_i|\leq \bar{U}$ for presentation simplicity in this paper. However, the presented strategies can be easily adapted to accommodate the case in which $-\underline{U}_i\leq u_i\leq \bar{U}_i$, where $\underline{U}_i$ and $\bar{U}_i$ are positive constants.
\end{Remark}

\begin{Remark}
Different from \cite{Mengscl13}-\cite{QiuIJRNC} that considered (optimal) consensus of multi-agent systems with bounded controls, this paper accommodates distributed Nash equilibrium seeking problems in systems with bounded controls. Compared with \cite{Mengscl13}\cite{YangAT14}, the problem is challenging as not only consensus of the players' estimates but also the optimization of the players' objective functions need to be achieved. In addition, the considered problem is challenging compared with \cite{XieSCL17} especially for second-order systems as it is difficult to distributively approximate $H(\mathbf{x})$ in \eqref{st_1}. Furthermore, \cite{QiuIJRNC} provided a projection operator based method to deal with distributed optimization problems in discrete-time systems with bounded controls, and hence, the design and analyses therein are distinct from this paper.
\end{Remark}

\section{Simulation studies}\label{p_1_numer}
This section verifies the effectiveness of the proposed seeking strategies in a mobile sensor network in which $x_i\in \mathbb{R}^2$ (denoted as $x_{i1}$ and $x_{i2}$, respectively). More specifically,  we consider the \textbf{connectivity control} for a network of  $3$ mobile sensors in which the sensors' objective functions are given by \cite{StankovicTAC12} $f_i(x_i,\mathbf{x}_{-i})=x_i^Tr_{ii}x_i+x_i^Tp_i+q_i+\sum_{j\in\mathcal{N}_i}m_{ij}||x_i-x_j||^2,$
where $r_{ii}\in \mathbb{R}^{2\times 2},p_i\in \mathbb{R}^{2\times 1},q_i\in \mathbb{R},m_{ij}\in \mathbb{R}$ are constant matrices, vectors or parameters and $\mathcal{N}_i$ denotes the physical neighboring set of player $i$. In the subsequent simulations, we consider Example 1 of \cite{StankovicTAC12} in which $i=3$, $r_{ii}$ for $i\in\{1,2,3\}$ are identity matrices,  and $m_{ij}=1$ except that $m_{13}=m_{31}=0.$ Moreover, $p_1=[2,-2]^T,p_2=[-2,-2]^T,p_3=[-4,2]^T,$ $q_i=3$ for $i\in\{1,2\}$ and $q_3=6.$ Through direct calculation, it can be easily verified that the example satisfies Assumptions \ref{Assu_1}, \ref{Assu_2}-\ref{Ass4} and the game admits a unique Nash equilibrium at $\mathbf{x}^*=[-0.125,0.75,0.75,0.5,1.375,-0.25]^T$ \cite{StankovicTAC12}.

In the following, velocity-actuated vehicles and acceleration-actuated vehicles will be simulated, successively.

\subsection{Velocity-actuated vehicles}
In this section, we consider velocity-actuated vehicles, whose dynamics can be described as $\dot{x}_i=u_i,$
where $x_i=[x_{i1},x_{i2}]^T$ denotes the position of sensor $i,$ $u_i=[u_{i1},u_{i2}]^T \in \mathbb{R}^2$, $u_{ij}$ for $i\in\{1,2,3\},j\in\{1,2\}$ denotes the control input of sensor $i$ that satisfies $|u_{ij}|\leq \bar{U}$.

\subsubsection{Saturated gradient play}\label{satu_ga}
\begin{figure}[t]
\begin{center}
\scalebox{0.38}{\includegraphics[86,383][554,491]{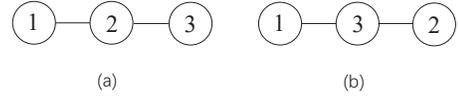}}
\caption{The communication graph among the sensors.}\label{comm_fig}
\end{center}
\end{figure}

\begin{figure}[htb!]
\begin{center}
\scalebox{0.38}{\includegraphics[10,229][578,608]{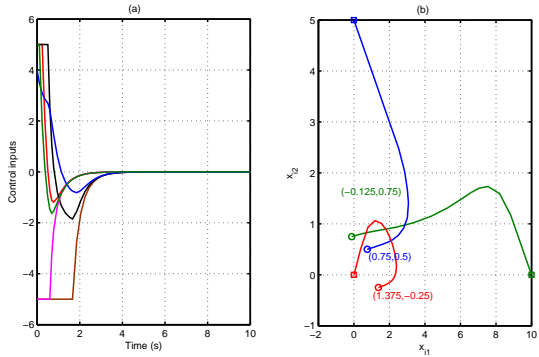}}
\caption{(a) and (b) show the control inputs and the trajectories of the sensors' positions   generated by the saturated gradient play in \eqref{satur_gra}, respectively.}\label{first_unb}
\end{center}
\end{figure}

In this section, we suppose that the mobile sensors can communicate with each other via the communication graph depicted in Fig. \ref{comm_fig} (a).
With $\mathbf{x}(0)=[10,0,0,5,0,0]^T$ and $\bar{U}=5,$ the trajectories of the sensors' positions and the control inputs generated by the saturated gradient play in \eqref{satur_gra} are depicted in Fig. \ref{first_unb}. Fig. \ref{first_unb} (a) illustrates that the control inputs are bounded by the given value and  Fig. \ref{first_unb} (b) shows that the sensors' positions would converge to the Nash equilibrium of the game asymptotically. Hence, by the simulation results, Theorem \ref{Lem_1} is numerically verified.

\subsubsection{Consensus-based distributed Nash equilibrium seeking}
In Section \ref{satu_ga}, the physical interactions among the sensors' objective functions coincide with their interactions in the communication graph. However, if this is not the case, the saturated gradient play can not be directly utilized in the distributed sensor networks. As an alternative, the distributed seeking strategy given in \eqref{bound_dis} can be adopted. To illustrate this case, in this section we suppose that the sensors can communicate with each other via the communication graph depicted in Fig. \ref{comm_fig} (b), which satisfies Assumption \ref{a3} as it is undirected and connected.

\begin{figure}[htb!]
\begin{center}
\scalebox{0.5}{\includegraphics[91,271][488,562]{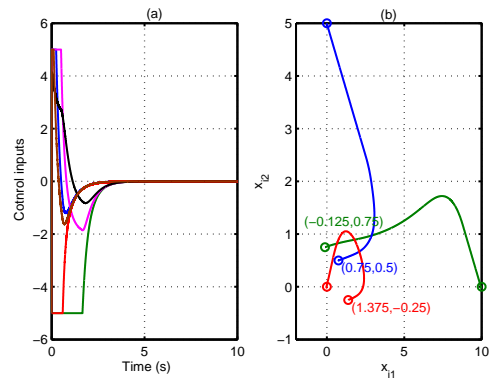}}
\caption{(a) and (b) show the control inputs and the trajectories of the sensors' positions  generated by the method in \eqref{bound_dis}, respectively.}\label{first_b}
\end{center}
\end{figure}

Let $\mathbf{x}(0)=[10,0,0,5,0,0]^T,$ $\bar{U}=5,$ $y_{ij}(0)=10$ and $\theta_{ij}=1000.$ By choosing $\mathcal{Q}$ and $\bar{\Theta}$ to be identity matrices, it can be verified that $\theta>\theta^*$, where $\theta^*$ is quantified in \eqref{lo_b}. Driven by the method in \eqref{bound_dis}, the control inputs are illustrated in Fig. \ref{first_b} (a) and the trajectories of the sensors' positions are plotted in Fig. \ref{first_b} (b). The control inputs stay within the bounded region as shown in Fig. \ref{first_b} (a). Moreover, Fig. \ref{first_b} (b) demonstrates that the trajectories of the sensors' positions would converge to the Nash equilibrium.  Hence, the effectiveness of  the proposed method in \eqref{bound_dis} is numerically verified.

%

\subsection{Acceleration-actuated vehicles}
\begin{figure}[htb!]
\begin{center}
\vspace{5mm}
\scalebox{0.38}{\includegraphics[98,272][491,566]{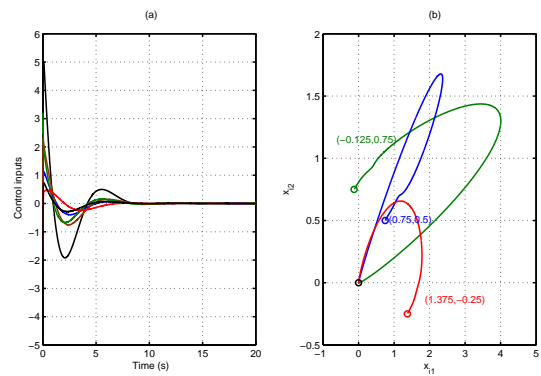}}
\vspace{5mm}
\caption{(a) and (b) show the control inputs and the trajectories of the sensors' positions  generated by the method in \eqref{method_4}, respectively.}\label{second_d}
\end{center}
\end{figure}

In this section, we suppose that the agents are acceleration-actuated vehicles whose dynamics can be described by $\dot{x}_i=\nu_i, \dot{\nu}_i=u_i,$
where $x_{i}=[x_{i1},x_{i2}]^T\in \mathbb{R}^2$ is the vector containing the positions of sensor $i,$ $\nu_{i}=[\nu_{i1},\nu_{i2}]^T\in \mathbb{R}^2$ is the vector containing the velocities of sensor $i$ and $u_i=[u_{i1},u_{i2}]^T\in\mathbb{R}^2$ is the vector containing the control inputs that satisfy $|u_{ij}|\leq \bar{U},$ for all $i\in\{1,2,3\},j\in\{1,2\}.$
Moreover, we suppose that the sensors update their positions according to \eqref{method_4}, in which $\bar{U}=5$, and all the variables are initialized at zero. Note that in the simulation, $\bar{K}_i=0.1$ and $\theta_{ij}=200$. By choosing $\mathcal{Q}$, $\bar{\Theta}$  to be  identity matrices and $K_i=0.1$, it can be verified that $\theta=200>\theta^*,$ where $\theta^*$ is defined in \eqref{b3}.  Under the communication graph depicted in Fig. \ref{comm_fig} (b), the simulation results are given in Fig. \ref{second_d}. As plotted in Fig. \ref{second_d} (a), the control inputs are bounded by the given value. Moreover, Fig. \ref{second_d} (b) depicts the trajectories of the sensors' positions, which shows that the sensors' positions asymptotically converge to the Nash equilibrium of the game.
The simulation results show that the proposed method in \eqref{method_4} is effective to achieve distributed Nash equilibrium seeking for second-order systems with bounded controls.

\section{Conclusions}\label{conc}
This paper considers Nash equilibrium seeking for games in systems where the control inputs are bounded. More specifically, first-order integrator-type systems are first considered, followed by second-order integrator-type systems. For both situations, we  first design a centralized seeking strategy based on the gradient play, which is further adapted to distributed networks. Based on the Lyapunov stability analysis, the convergence properties of the designed algorithms are analytically investigated. It is shown that the proposed seeking strategies would enable the players' actions to converge to the Nash equilibrium under the given conditions.

\end{document}